\documentclass[12pt,reqno]{amsart}
\usepackage{pgfplots}
\usepackage{tikz}
\pgfplotsset{compat=newest}
\pgfplotsset{ytick style={draw=none}}
\pgfplotsset{xtick style={draw=none}}
\usepackage{amsmath,amsfonts,amssymb,amsthm}
\usepackage{graphicx}
\graphicspath{ }
\usepackage{mathrsfs}
\usepackage{epstopdf}
\usepackage{csquotes}
\usepackage{wrapfig}
\usepackage{accents}
\usepackage{caption}
\usepackage{subcaption}
\usepackage{calligra}
\usepackage[colorlinks]{hyperref}
\usepackage{kantlipsum}
\usepackage{cleveref}
\hypersetup{colorlinks, 
	breaklinks,
	linkcolor=red,
	citecolor=red,
	linktocpage=true}

\numberwithin{figure}{section}

\theoremstyle{plain}
\newtheorem{thm}{Theorem}[section]

\newtheorem{prop}[thm]{Proposition}

\theoremstyle{definition}

\theoremstyle{remark}

\usepackage{mathtools}
\makeatletter
\@namedef{subjclassname@2020}{%
	\textup{2020} Mathematics Subject Classification}
\makeatother
\title[FUCHSIAN SCHOTTKY GROUP]{CONSTRUCTION OF FUCHSIAN SCHOTTKY GROUP WITH CONFORMAL BOUNDARY AT INFINITY }

\author[A. A. Shaikh]{Absos Ali Shaikh$^{*1}$}
\author[U. Roy]{Uddhab Roy$^{2}$}

\address{$^{1,2}$Department of Mathematics\\ The University of Burdwan\\ Burdwan--713104\\ West Bengal\\ India.}

\email{$^1$aask2003@yahoo.co.in, aashaikh@math.buruniv.ac.in}
\email{$^2$uddhabroy2018@gmail.com}

\begin{document}

\noindent\footnotetext{$^*$ Corresponding author.\\
	$\mathbf{2020}$  \textit{Mathematics Subject Classification}.  Primary 20H10; Secondary 57K20, 30F35, 30F40, 37F32.}

\maketitle 	

%\begin{center}
 %        (Dedicated to the memory of Professor Lajos Tam\'assy)

%\end{center}

\begin{abstract}
In this article, we have constructed an interesting type of generalized Schottky group, named as Fuchsian Schottky group of arbitrary finite rank, in the context of the classical Schottky group (i.e., Schottky curves which are Euclidean circles). After that, we initiated the construction of the generalized Fuchsian Schottky group of any finite rank by including orientation-reversing isometries of the hyperbolic plane as side-pairing transformations. Further, we have investigated the hyperbolic ends for any arbitrary finite rank Fuchsian Schottky groups from the point of view of the Euler characteristic in the hyperbolic surface. Finally, we have shown that the compact core of the conformally compact Riemann surface can be decomposed into non-tight pairs of pants by using suitable twist parameters with some fixed Bers' constant. The Fenchel-Nielsen coordinates for Teichm\"uller space corresponding to any finite rank Fuchsian Schottky groups are also obtained.

\end{abstract}

\section{\textbf{INTRODUCTION}}

 It is well known that Schottky groups are free, discrete subgroups of $PSL(2, \mathbb{C})$. Since the latter half of the $20$th century, Schottky groups have been studied by various authors like Chuckrow (\cite{Chuckrow} and \cite{Chuckrow1970}), Marden \cite{Marden}, Ber \cite{Ber}, Button \cite{Button}, Peter \cite{Peter}, and Maskit (\cite{Maskit}, \cite{Maskit1967}, and \cite{Maskit2000}). The classical construction of this group is as follows: consider the region in $\mathbb{CP}^1$ bounded by $n$ pairs of mutually disjoint circles and pair these sides by $n$ loxodromic transformations such that the quotient is a hyperbolic surface of genus $n$; the Schottky group is the free group generated by these side-pairing transformations.
  A similar construction also works in $PSL(2, \mathbb{R})$, taking a region D of the hyperbolic plane bounded by semi-circles orthogonal to the boundary at infinity; these are called real Schottky groups. Now, to develop the literature, in this manuscript, we have constructed a type of generalized Schottky group (named as Fuchsian Schottky group) on the analogy of the real Schottky groups with two additional conditions:
  
  (i) The transformation is paired with the semi-circles with its reflection on the upper imaginary axis.
  
  (ii) The positions of semi-circles at the circle at infinity are non-tangential (at least for the Fuchsian Schottky group). \\ After that, we built the construction of the generalized Fuchsian Schottky groups of finite rank by including the orientation-reversing isometries as side-pairing transformations in the upper-half plane model. Then we have derived the hyperbolic ends for any finite rank Fuchsian Schottky groups from the point of view of the topological invariant Euler characteristic in the hyperbolic surface. This study gives rise to interesting surfaces, like finite Loch Ness monster and finite Jacob's ladder with sufficiently small limit sets. In fact, in this paper, we fabricated a finitely generated hyperbolic group which produces an infinite area hyperbolic surface. Further, we have deduced the Fenchel-Nielsen coordinates for Teichm\"uller space corresponding to any arbitrary finite rank Fuchsian Schottky groups by indicating the Bers' constant.\\

 \textbf{Outline of this paper.} In section 2 we have briefly discussed the preliminaries of the generalized Schottky group and the geometry of ends for hyperbolic surfaces. In section 3 we have constructed a rank $n$, $n \in$ $\mathbb{N} - \{1\}$, purely hyperbolic generalized Schottky group, called the Fuchsian Schottky group, (denoted by $\Gamma_{S_n}$) by using orientation-preserving isometries of $\mathbb{H}$ as side-pairing transformations (see, Proposition $3.1$). In section 4 we have extended this construction in $PSL^*(2,\mathbb{R})$ (see, $Table:1$, for description) by including orientation-reversing isometries of $\mathbb{H}$ as side-pairing transformations (see, Proposition $4.1$). In this paper, we have used orientation-reversing isometries only in section 4. Then we studied the geometry of the limit set for the group $\Gamma_{S_n}$, $n \in \mathbb{N}$. In subsection  5.1, we have organized surface decomposition and gluing for $\Gamma_{S_n}$, $n \in \mathbb{N}$, including the discussion about the equivalency of the convex core and the compact core for the group $\Gamma_{S_n}$. Then we have given a characterization of the Nielsen region for any finite rank Fuchsian Schottky group, in the hyperbolic plane. In subsection 5.2, we have investigated the quotient surfaces with the conformal boundary at infinity corresponding to $\Gamma_{S_n}$ ($n \in \mathbb{N}$) from the point of view of the Euler characteristic, $\chi$, in the hyperbolic surface (see, Theorem $5.1$ and Theorem $5.2$). In section 6 we have first discussed the existence of non-tight pairs of pants for $\Gamma_{S_n}$, $n \in \mathbb{N} - \{1\}$. After that, we created hyperbolic surfaces attaching $Y$ and $X$-pieces by using $(3n-2)$ numbers of twist parameters for rank $n$ Fuchsian Schottky group, $n \in \mathbb{N} - \{1\}$. A characterization of half-collars for $\Gamma_{S_n}$ is also studied. Finally, we have proved that the compact core of the conformally compact Riemann surface can be decomposed into $(2n-2)$ numbers of non-tight pair of pants where Bers' constant is strictly less than $(31n+21)$ and produced Fenchel-Nielsen coordinates $\{l_1, l_2, ..., l_{(3n+2)}$ ; $\beta_1, \beta_2, ..., \beta_{(3n-2)}\}$ of Teichm\"uller space ($\subset$ $\mathbb{R}^{6n-4}$) corresponding to the Fuchsian Schottky group of rank $n$, $n \in \mathbb{N} - \{1\}$ (see, Theorem $6.1$).

  \section{\textbf{GENERALIZED SCHOTTKY GROUP AND ENDS}}
 
   In this section, we have discussed some basics of the generalized Schottky group, hyperbolic ends, Y-piece, and conformal boundary for the non-abelian Fuchsian group. For details, we refer the readers to the books of Beardon \cite{Beardon}, Maskit \cite{Maskit1988}, Matsuzaki \cite{Matsuzaki1998}, and Buser \cite{Buser}.
   
    A generalized Schottky group $S(\gamma_1,\gamma_2, ..., \gamma_p)$ (with respect to the point $0$) is a subgroup of $PSL(2,\mathbb{R})$ generated by non-trivial non-elliptic elements \{$\gamma_1, \gamma_2, ..., \gamma_q$\}, $(1 \le p \le q)$ such that for all $m \ne n$ we have $$\{(C(\gamma_m) \cup C(\gamma^{-1}_m)) \cap (C(\gamma_n) \cup C(\gamma^{-1}_n))\} = \phi,$$ where $C(\gamma_i) = \{a\in \mathbb{H} : d(a, \gamma(0)) \le d(a, 0)\} $, $d$ denotes the hyperbolic distance. The sets $(C(\gamma_m) \cup C(\gamma^{-1}_m))$ are pairwise tangential if $\gamma_m$ are parabolic, however, these sets are pairwise non-tangential when $\gamma_m$ are hyperbolic. Notice that, if we omit the tangential case then the group $S(\gamma_1,\gamma_2, ..., \gamma_p)$ is reduced to a Schottky group in $PSL(2,\mathbb{R})$. Although it is still possible that $C(\gamma_m)$ and $C(\gamma^{-1}_m)$ are tangent to each other. For the generalized Schottky group, the parabolic isometry is not always conjugate to powers of the parabolic generator $\gamma_i$, whereas, for the Schottky group, it always agrees.
   
   The following $Table:1$ represents the symbols and notations of different terms that will be used throughout the paper.
                \begin{center}
                \resizebox{16cm}{4.1cm}{\begin{tabular}{|| c c c c ||} 
                                                  \hline
                                                  Sl. No. & Descriptions of terms & & Notation \\ [0.5ex] 
                                                  \hline \hline
                                                  1. & hyperbolic plane &  & $\mathbb{H} $\\ 
                                                  \hline
                                                  2. & convex core &  & $C$ \\
                                                  \hline
                                                  3. & limit set &  & $\Lambda(\Gamma)$ \\
                                                  \hline
                                                  4. & set of discontinuity &  & $\Omega(\Gamma)$ \\
                                                  \hline
                                                  5. & projective special linear group over $\mathbb{R}$ degree $2$ &  & $PSL(2,\mathbb{R}$)\\
                                                  \hline
                                                  
                                                  6. & boundary of the hyperbolic plane &  &  $\partial\mathbb{H}$\\
                                                  \hline
                                                  7. & Nielsen region &  & $N$ \\
                                                  \hline
                                                  8. & compact core &  & $K$ \\
                                                  \hline
                                                  
                                                  9. & Euler's characteristic & & $\chi$ \\
                                                  \hline
                                                  10. & $SL^{*}(2, \mathbb{R})$ is a group of real $2\times 2$ matrices with 
                                                  
                                                   $det=$ $\pm 1$, $SL^{*}(2, \mathbb{R})/\{\pm I_2\}=$ & &  $PSL^{*}(2,\mathbb R)$   \\
                                                  \hline
                                                  11. & circle at infinity & & $\mathbb{RP}^1$ \\
                                                  \hline
                                                  12. & Fuchsian Schottky group of rank $n$ & & $\Gamma_{S_n}$ \\
                                                  \hline
                                                  
                                                  13. & conical limit set &  &  $\Lambda_c(\Gamma)$ \\
                                                  \hline
                                                  14. & Hausdorff dimension &  &  $dim_H$  \\
                                                  \hline
                                                  15. & empty set &  &  $\phi$ \\
                                                  \hline
                                                  16. & $\mathbb{H} \cup \{\infty\}$ &  &  $\overline{\mathbb{H}}$ \\
                                                  \hline
                                                  17.  & the conformal boundary at infinity &  & $\Omega(\Gamma)/\Gamma$ \\
                                                                               \hline
                                                  18. & funnels end & & $F'$ \\ 
                                                                            \hline
                                                                            19. & Sphere at infinity &  & $\mathbb{CP}^1$ \\
                                                                            \hline
                                                    20. & projective special linear group over $\mathbb{C}$ degree $2$ & & $PSL(2,\mathbb{C}$) \\
                                                                                                                                         \hline                     
                                                 \end{tabular}}

                                 \end{center}
                                 
                                 $$ Table: 1$$

    All surfaces in
    this paper are orientable and triangulable. A funnel is a hyperbolic surface having one geodesic boundary
   component which is isometric to $ D/ < z \rightarrow e^{l_1}z >$, where $ D = \{z \in U : \Re(z) \le 0 \} $ has
   the induced metric as a subspace of the upper-half plane model U of the hyperbolic plane ($\Re(z)$ denotes the real part of $z=x+iy$). A cusp is a hyperbolic surface having one horocyclic boundary component
   which is isometric to the quotient $\{z: \Im(z) \ge 1 \}/ < z \rightarrow z +1 > $, where  $\{ z: \Im(z)\ge 1 \}$
    has the induced metric as a subspace of the upper-half plane model ($\Im(z)$ denotes the imaginary part of $z$). We all know that any finite-type geodesically complete hyperbolic surface has ends that are either funnels or cusps. A geodesically complete hyperbolic surface $ X $ is the quotient
    of the hyperbolic plane $\mathbb{H}$, by a  non-elementary torsion-free Fuchsian  group $\Gamma$. Throughout this paper, we have used the Fuchsian group to mean the torsion-free Fuchsian group. The action of $\Gamma$ on the circle at infinity of the hyperbolic
    plane breaks up into the limit set $\Lambda(\Gamma)$ and the (if not empty) set of discontinuity, $\Omega(\Gamma)$. The set of discontinuity $\Omega(\Gamma)$ is made up of a countable union of intervals of discontinuity. The stabilizer in the finitely generated Fuchsian group $\Gamma$  of an interval of discontinuity is generated
    by a hyperbolic element. The convex core of $ X$, i.e., $C(X)$, is the quotient of the convex hull of the limit set. The convex core is the smallest closed convex sub-surface with the boundary
    which carries all the homotopy. Let, $X = \mathbb{H}/\Gamma$ be a geodesically complete hyperbolic surface. $\Gamma$ is said to be of the first kind if $\Lambda(\Gamma) = \partial\mathbb{H}$. Otherwise it is of the second kind. Similarly, $X$ is said to be of the first kind when $C(X) = X$. Otherwise it is of the second kind that we are mainly interested in this manuscript. An end is usually called geometrically finite if it contains a neighborhood that is disjoint from $N/\Gamma$ and geometrically infinite otherwise. If  $X$ is a geodesically complete
    hyperbolic surface, then a boundary component of $C(X)$ is either a simple closed geodesic that bounds a funnel in $X$ and corresponds to an isolated
    end that is visible or a simple geodesic isometric to the real line that bounds a half-plane in $X$ and
    corresponds to a component of a visible end of infinite type. Also, a visible infinite type end of a geodesically complete hyperbolic
    surface  has an equivalence class of components of the convex core boundary being simple
    open geodesics with attached half-planes. A surface may have an uncountable number of ends. However, the hyperbolic metric in the hyperbolic plane 
    places restrictions on the geometry of the ends. Since a Fuchsian group has
    only a countable number of intervals of discontinuity a complete hyperbolic
    surface has at most a countable set of ends (that means at most enumerable) that are visible (see, \cite{Basmajian2019} for more details). A Fuchsian group $\Gamma$ is said to be a pair of pants if the quotient space $\mathbb{H}/\Gamma$ is topologically a sphere with $3$ disks removed. A pair of pants (or a $Y$-piece) is a compact Riemann surface of signature ($0,3$) whereas an $X$-piece is a compact Riemann surface of signature ($0,4$). A non-tight hyperbolic pair of pants is a hyperbolic sphere with three geometric conformal holes, where all the conformal holes are geodesic boundary components as simple closed curves (i.e., homeomorphic to three circles). A compact Riemann surface is a hyperbolic surface that is compact without a boundary. A compact hyperbolic surface of signature ($g,m$) is called a Riemann surface of signature ($g,m$) when every boundary component is a smooth closed geodesic. For a torsion-free non-elementary Fuchsian group $\Gamma$, the conformal boundary of a hyperbolic $2$-manifold is the topological boundary of $(\mathbb{H}\cup \Omega(\Gamma))/\Gamma$. The conformal boundary at infinity of the hyperbolic surface is defined by $\Omega(\Gamma)/\Gamma$, for the non-abelian group $\Gamma$, where this group acts as a group of isometries. Canary \cite{Canary} has proved that the length of a curve in the conformal boundary produces an upper bound on the length of the corresponding curve in the convex core boundary.

\section{\textbf{CONSTRUCTION OF ARBITRARY FINITE RANK FUCHSIAN SCHOTTKY GROUP IN THE CONTEXT OF THE CLASSICAL SCHOTTKY GROUP}}
  In this section we have constructed a purely hyperbolic generalized Schottky group, called the Fuchsian Schottky group of rank $n$, $n$ $\in$ $\mathbb{N} - \{1\}$, using $2n$ semi-circles in the hyperbolic plane, $\mathbb{H}$, with centers lying on the real projective line, $\mathbb{RP}^1$.\\ 
    \begin{prop} \label{p3.1}
     Any rank $n$ $\in$ $\mathbb{N} - \{1\}$, Fuchsian Schottky group contains orientation preserving isometries of $\mathbb{H}$ as  side-pairing transformations.
    \end{prop}
   
   \begin{proof}
   In the upper-half plane model, let $D_{1},D_{2},...,D_{2n-1}$, and $D_{2n}$ be open mutually disjoint Euclidean semi-circles in $\mathbb{H}$ with centers on the real projective line, $\mathbb{RP}^1$. The transformation is paired with the semi-circle with the diameter on the real axis and with its reflection on the upper imaginary axis. For each pair, $D_j$, $D_{j+n}$, $(j=1,2,...,n)$ we have supposed $S'_j$ $\in$ $PSL(2, \mathbb{R}$) is a transformation that sends $\partial$$D_j$ to $\partial$$D_{j+n}$ ($\partial$$D_j$ denotes the boundary of $D_j$). Also, the transformation $S'_{j}$ maps the (whole) exterior of $D_{j}$  to the interior of $D_{j+n}$. Each transformation $S'_{j}$ is hyperbolic with a repelling fixed point inside $D_{j}$ and an attracting fixed point inside $D_{j+n}$.
  
      Let, $D$ be a convex non-compact hyperbolic polygon with free edges (a polygon has edges that are lying on the boundary, such edges are called free edges).  $D$ has vertices at the  boundary at infinity $\{a_{1}, a_{2}, ..., a_{n}$; $-a_{1}, -a_{2}, ..., -a_{n}$; $0$; $b_{1}, b_{2}, ..., b_{n}$; $-b_{1}, -b_{2}, ..., -b_{n}\}$. Notice that, here every vertex of the polygon $D$ is an improper vertex except $0$, which is a neither proper nor improper vertex (see, \cite{Beardon}). The sides of $D$ are the geodesic segments joining $a_{i}$ and $b_{i}$; $-b_{i}$ and $-a_{i}$. Assume that all sides $s_{i}$ and $s'_{i}$, $\{i=1,2,...,(n-1),n\}$ of $D$ are equipped with side-pairing transformations $\gamma_{i}$ and $\gamma'_{i}$ respectively. We consider the orientation-preserving isometry of $\mathbb{H}$ to be a side-pairing transformation. The isometries $\{ \gamma_{1},\gamma_1',\gamma_{2},\gamma_{2}',...,\gamma_{n-1},\gamma_{n-1}',\gamma_{n},\gamma_{n}' \}$ (which are orientation preserving) act locally such that the half-plane bounded by $\{s_{1}, s_{2}, ..., s_{n-1},$
       $s_{n}$ ; $s'_{1}, s'_{2}, ..., s'_{n-1},$ $ s'_{n} \}$ containing $D$ is mapped by $\gamma_{i}$ and $\gamma'_{i}$ to the half-plane bounded by $\gamma_{i}(s_i)$ and $\gamma'_{i}(s'_i)$ respectively. Let, each free edge be paired with itself by the identity map. Now, given a vertex (let, $a_{i}$) and a side (let, $s_{i}$ or $s'_{i}$) with an endpoint at that vertex, we have the following cases:

               (a) (i) The vertex $a_{i}$ belongs to a parabolic cycle which does not contain a vertex that is the end-point of a free edge. (ii) The vertex $a_{i}$ $\in$ $\partial\mathbb{H}$ belongs to a free cycle if the parabolic cycle does contain a vertex that is the end-point of a free edge.
                   
               (b) The vertex $a_{i}$ $\in$ $\partial\mathbb{H}$ belongs to a free cycle if the hyperbolic cycle does contain a vertex that is the end-point of a free edge. \\ Note that, no vertex can belong to the parabolic cycle.
         
               Now, to show the discreteness of $\Gamma_{S_n}$ (say) in $PSL(2,\mathbb{R})$ we have used Poincar\'e's theorem (see, \cite{Beardon}). Observe that, it is not usual to apply when all the edges of the polygon are free. But we have presented that as follows: let, $D$ be a convex hyperbolic polygon with all edges free. Suppose that $D$ $\subset$ $\mathbb H$ is a finite-sided convex polygon whose sides are identified in pairings by isometries  $\{s_{1}, s_{2}, ..., s_{n-1}, s_{n}$ ; $s'_{1}, s'_{2}, ..., s'_{n-1}, s'_{n} \}$ = $G_{1}$. Now, $D$ is equipped with a collection $G_{1}$ of side-pairing hyperbolic orientation-preserving isometries. Assume that, each free edge is paired with itself via the identity. Also, each hyperbolic cycle satisfies the hyperbolic cycle condition. Then, the subgroup $<G_{1}>$ generated by $G_{1}$ is a subgroup of $PSL(2,\mathbb R)$ and $D$ is a fundamental domain of $G_{1}$. Now, it is clear that all vertices lie on the boundary and every vertex is the endpoint of a free edge. So, as every vertex must belong to a  free cycle, the side-pairing transformations generate a discrete subgroup ($<G_{1}>$ $=$ $\Gamma_{S_n}$, say) of $PSL(2,\mathbb R)$. Here, $\gamma_{1},\gamma_1',\gamma_{2},\gamma_{2}',...,\gamma_{n-1},\gamma_{n-1}',\gamma_{n},\gamma_{n}'$ are orientation preserving isometries given by the hyperbolic transformations with fixed points at $a_{1},b_{1}$ ; $-a_{1},-b_{1}$ ; $a_{2},b_{2}$ ; $-a_{2},-b_{2}$ ; ... ; $ a_{n-1},b_{n-1}$ ; $-a_{n-1},-b_{n-1}$ ; $ a_{n},b_{n}$ ;$-a_{n}$,$-b_{n}$ respectively $\in\partial\mathbb{H}$. Clearly, $\gamma_{1},\gamma_{1}',\gamma_{2},\gamma_{2}', ...,
              $ $ \gamma_{n-1},\gamma_{n-1}',\gamma_{n},\gamma_{n}'$ are side-pairing transformations. Now, it is obvious that cycles containing the vertices are free cycles and there is no parabolic cycle containing any of the above vertices. So, all the cycles are hyperbolic cycles and all transformations are hyperbolic transformations. Hence, by {Poincar\'e}'s theorem, $\{\gamma_{1},\gamma'_{1},\gamma_{2},\gamma'_{2}, ..., \gamma_{n-1},\gamma'_{n-1},\gamma_{n},\gamma'_{n}\}$ generates a discrete subgroup of $PSL (2,\mathbb R)$. Observe that, for the group     $<\gamma_{1},\gamma'_{1},\gamma_{2},\gamma'_{2}, ..., \gamma_{n-1},\gamma'_{n-1},\gamma_{n},\gamma'_{n}>$, the generators are overlapping by $\gamma_{1},\gamma'_{1}$ ; $ \gamma_{2},\gamma'_{2}$ ; $ ... $ ; $ \gamma_{n-1},\gamma'_{n-1}$ ; $\gamma_{n},\gamma'_{n}$. Hence, this group reduces to the group $<\gamma_{1},\gamma_{2}, ..., \gamma_{n-1}, \gamma_{n}>$ = $\Gamma_{S_n}$ (say). This group is a Fuchsian group with $n$ generating elements. Here, $D$ is a fundamental domain for $G_{1}$ and $D$ is of infinite area. Now, it is well known that (see, \cite{Beardon}) a Fuchsian group is of the first kind if and only if a (hence all) fundamental domain has a finite area. So, $\Gamma_{S_n}$ must be a Fuchsian group of the second kind. Hence, the limit set of $\Gamma_{S_n}$, i.e., $\Lambda(\Gamma_{S_n})$ contains either $0,1,$ or $2$ elements or a Cantor set. Now, $\gamma_{1}$ has hyperbolic fixed points at $a_{1}$ and $b_{1}$. Hence, $a_{1}$, $b_{1}$ $\in$  $\Lambda(\Gamma_{S_n})$. Similarly, $\Gamma_{S_n}$ contains the hyperbolic M\"obius transformations $\{\gamma_{2},\gamma_{3}, ..., \gamma_{n-1},\gamma_{n}\}$ having fixed points at $\{a_{2},b_{2}$; $a_{3},b_{3}$; ... ;$a_{n-1},b_{n-1}$; $a_{n},b_{n}\}$ respectively. 
              So, $a_1,b_1,a_2,b_2, ..., a_{n-1},b_{n-1},a_n,b_n$ $\in$ $\Lambda(\Gamma_{S_n})$.
              Similarly, one can choose $\{\gamma'_1, \gamma'_2, ...,  \gamma'_{n-1}, \gamma'_{n}\}$ as the generators of the group $\Gamma_{S_n}$. Then for that case, $-b_1, -a_1, -b_2, -a_2, ..., -b_{n-1}, -a_{n-1},
              -b_n, -a_n$ $\in$ $\Lambda(\Gamma_{S_n})$. Therefore, $\Lambda(\Gamma_{S_n})$ contains more than 2 elements. That means $\Lambda(\Gamma_{S_n})$ is a Cantor set. We call the group $\Gamma_{S_n}$ is the rank $n$ Fuchsian Schottky group. This completes the proof.
               \end{proof}
                
              In this way, one can easily construct any rank $n$, $n$ $\in \mathbb{N} - \{1\}$, Fuchsian Schottky group by introducing  Euclidean open mutually disjoint semi-circles $D_i$, \{$i=1,2,..., (2n-1), 2n$\}  in $\mathbb{H}$ centers on $\mathbb{RP}^1$ with the number of semi-circles are double of the required number of rank of that group. 
              Also, for $\Gamma_{S_n}$, each $j$ $\in$ $\{1,2, ... ,n\}$,
              $S'_{j+2n}$ $=$ $S'_j$  and  $S'_{j+n}$ $=$ $(S'_j)^{-1}$ hold, where the elements $\gamma_s$ or $\gamma'_s$ ($\in$ $\Gamma_{S_n}$) are given by fractional linear transformations that map the whole exterior of $D_k$ to the interior of $D_t$ with $\lvert k-t \rvert$ $=$ $n$, where `$n$' is the rank of the group $\Gamma_{S_n}$.\\
             
              \textbf{Remark 1.}
               Let, $X^\circ$ be the interior of a connected orientable hyperbolic surface with a boundary. Then the surface group of $X^\circ$, i.e., the fundamental group $\pi_1(X^\circ)$ is free on ($1 - \chi$) generators for $\Gamma_{S_n}$. The above constructed group $\Gamma_{S_n}$ ($\subset$ $PSL(2,\mathbb{R})$) is the holonomy of an infinite area hyperbolization of $X^\circ$. So, our Fuchsian Schottky representation is not exhaustive. Here, $\Gamma_{S_n}$ is freely generated by $S'_1, S'_2, ..., S'_n$ $\in$ $PSL(2, \mathbb{R})$ and there exists $2n$ open intervals \{$(a_1, b_1), (-b_1, -a_1)$  ; $(a_2, b_2), (-b_2, -a_2)$  ; ... ; $(a_n, b_n), (-b_n, -a_n)$\} which are  mutually disjoint  and centers on the real projective line such that $S'_k (-b_k, -a_k) = (a_k, b_k)$, with  \{$(\overline{I^+_1} \cup \overline{I^+_2} \cup ... \cup \overline{I^+_k})$ $\cup$ $(\overline{I^-_1} \cup \overline{I^-_2} \cup ... \cup \overline{I^-_k})$\}  $\subsetneq$ $\mathbb{RP}^1$, where $I^+_k = (a_k, b_k)$ and $I^-_k = (-b_k, -a_k)$. We observe that this group $\Gamma_{S_n}$, $n$ $\in$ $\mathbb{N} - \{1\}$, is not a lattice, although both $\Gamma_{S_n}$ and lattices are geometrically finite.

               \section{\textbf{CONSTRUCTION OF ARBITRARY FINITE RANK GENERALIZED FUCHSIAN SCHOTTKY GROUP IN THE CONTEXT OF THE CLASSICAL SCHOTTKY GROUP}} 
               
               In section $3$ we have utilized only orientation-preserving isometry to be a side-pairing transformation. But in this section, we have included the orientation reversing isometry of $\mathbb{H}$ to be a side-pairing transformation to extend the Proposition $3.1$ construction in $PSL^*(2,\mathbb{R})$. We have concluded this section by briefly describing the geometry of the limit set for this generalized Schottky group $\Gamma_{S_n}$, $n$ $\in$ $\mathbb{N}$.\\
                     
                       \begin{prop}\label{p4.1}  Any rank $n$, $n \in \mathbb{N} - \{1\}$, generalized Fuchsian Schottky group contains orientation preserving and orientation reversing isometries of $\mathbb{H}$ as  side-pairing transformations.
                       \end{prop}

                      \begin{proof}
                       At first, taking the notion of our Fuchsian Schottky construction in Proposition $3.1$, we have allowed arbitrary isometry of $\mathbb{H} $ to be a side-pairing transformation (i.e., orientation reversing isometry may occur). Now, let the sides of $D$ be the geodesic segments joining $a_{i}$ and $b_{i}$; $-a_{i}$ and $-b_{i}$ but all sides $t_{i}$ and $t'_{i}$, $\{i=1,2,...,(n-1),n\}$ of $D$ are equipped with side-pairing transformations $\gamma_{i}$ and $\gamma'_{i}$ respectively which may be orientation-preserving or orientation-reversing isometries. So, the isometries $\{ \gamma_{1},\gamma_1',\gamma_{2},\gamma_{2}',...,\gamma_{n-1},\gamma_{n-1}',\gamma_{n},\gamma_{n}' \}$ act locally such that the half-plane bounded by $\{t_{1}, t_{2}, ..., t_{n-1}, t_{n}$ ; $t'_{1}, t'_{2}, ..., t'_{n-1},$ $ t'_{n} \}$ containing $D$ is mapped by $\gamma_{i}$ and $\gamma'_{i}$ to the half-plane bounded by $\gamma_{i}(t_i)$ and $\gamma'_{i}(t_i)$ respectively but opposite $D$. Also let, each free edge be paired with itself by the identity map. Now, to show the discreteness of $<G_2>$ (say) in $PSL^{*}(2, \mathbb{R})$ we have again utilized Poincar\'e's theorem as follows: let, $D$ be a convex hyperbolic polygon consisting of all free edges. Assume that, $D$ $\subset$  $\mathbb{H}$ is a finite-sided convex polygon whose sides are identified in pairings by isometries $\{t_{1}, t_{2}, ..., t_{n-1}, t_{n}$; $t'_{1}, t'_{2}, ..., t'_{n-1},$ $ t'_{n} \}$ $=$ $G_2$. So, $D$ is equipped with a collection $G_2$ of side-pairing hyperbolic isometries where each free edge is paired with itself via the identity. Suppose that, each hyperbolic cycle satisfies the hyperbolic cycle condition. Then, the subgroup $<G_2>$ generated by $G_2$ is a discrete subgroup of $PSL^{*}(2,\mathbb R)$. Observe that, $D$ is not a fundamental domain of $G_2$, because we have taken arbitrary isometries (may not just M\"{o}bius transformations) as side-pairing transformations for $G_2$. In Proposition $3.1$, we have assumed that $\gamma_{1},\gamma_1',...,\gamma_{n},\gamma_{n}'$ are orientation-preserving isometries which are given by the hyperbolic transformations with fixed points at $a_{1},b_{1}$; $-b_{1},-a_{1}$; ...; $a_{n},b_{n}$; $-b_{n},-a_{n}$ respectively $\in\partial\mathbb{H}$. Now, let some of the side-pairing transformations be orientation-reversing isometries. Then, the group $\Gamma^1 $ $=$ \{$<G_2> \cap$ $ PSL(2,\mathbb R)\}$, of all orientation preserving transformations in $<G_2>$ is a Fuchsian group. Let, $\gamma_i$ $\in$ $<G_2>$ be an orientation reversing isometry, where $i = 1,2, ..., (n-1)$ or $n$. Then, $D_1$ $=$ \{$D$ $\cup$ $\gamma_i(D)\}$ is a fundamental domain for $\Gamma^1$. So, $D_1$ $\supset$ $D$. Notice that, all the vertices of the hyperbolic polygon lie on the boundary of the hyperbolic plane. Also, all vertices of the polygon are the endpoints of the free edges. Therefore, all vertices belong to free cycles. Hence, the side-pairing transformations generate a discrete subgroup $<G_2>$ = $\Gamma^2$ (say) of $PSL^*(2,\mathbb{R})$. Again, the orientation reversing isometry $\gamma_i$ given by the composition of a hyperbolic transformation with fixed points at $a_i$, $b_i$ $\in$ $\partial\mathbb{H}$ and a reflection in the geodesic $[a_i, b_i]$. Also, $\gamma_i$ pairs the side $[a_i, b_i]$ to itself. So, $\gamma_i$ is a side-pairing transformation. For orientation-preserving isometries, it is obvious that they are side-pairing transformations. Now, for the orientation reversing case, we have clarified this in the ensuing way. \\ The cycle containing the vertex $(a_{i-1})$, $i = \{2,3, ..., n$ or $ (n+1)\}$ is given by ;$$\binom{a_{i-1}}{t_{i-1}}\xrightarrow{\gamma_{(i-1)}}\binom{a_{i-1}}{t_{i-1}}  \xrightarrow{*}{\binom{a_{i-1}}{t_{(i-2)(i-1)}} } \xrightarrow{id}{ \binom{a_{i-1}}{t_{(i-2)(i-1)}}}  \xrightarrow{*} \binom{a_{i-1}}{t_{i-1}}.$$  As this cycle contains the free edge $t_{(i-2)(i-1)}$, it is a free cycle. \\ And the cycle containing the vertex $(a_i)$, $i=1$ (only) is given by ;
                               
                               $$\binom{a_{i}}{t_{i}}\xrightarrow{\gamma_{(i)}}\binom{a_{i}}{t_{i}}  \xrightarrow{*}{\binom{a_{i}}{t'_{(i)(i)}} } \xrightarrow{id}{ \binom{a_{i}}{t'_{(i)(i)}}}  \xrightarrow{*} \binom{a_{i}}{t_{i}}.$$  As this cycle contains the free edge $t'_{(i)(i)}$, it is a free cycle. \\ Now, the cycle containing the vertex $(-b_i)$, $i = \{1, 2, ..., n$ or $ (n-1)\}$ is given by ; $$\binom{-b_{i}}{t'_i}\xrightarrow{\gamma'_n}\binom{-b_{i}}{t'_i}  \xrightarrow{*}{\binom{-b_{i}}{t'_{(i+1)i}} } \xrightarrow{id}{ \binom{-b_{i}}{t'_{(i+1)i}}}  \xrightarrow{*} \binom{-b_{i}}{t'_{i}}.$$ \\ As this cycle contains the free edge $t'_{(i+1)i}$, it is also a free cycle. \\ And for the cycle containing the vertex $(-b_i)$, $i=n$ (only) is given by ;
                                $$\binom{-b_{i}}{t'_i}\xrightarrow{\gamma'_n}\binom{-b_{i}}{t'_i}  \xrightarrow{*}{\binom{-b_{i}}{t'_{(i)(i)}} } \xrightarrow{id}{ \binom{-b_{i}}{t'_{(i)(i)}}}  \xrightarrow{*} \binom{-b_{i}}{t'_{i}}.$$ \\ We notice that, for $i=n$, the free edge $t'_{(i)(i)}$ is equal to $t_{(i)(i)}$. Now, as this cycle contains the free edge $t'_{(i)(i)}$, it is also a free cycle. Clearly, there is no parabolic cycle. So, all the cycles are hyperbolic cycles and all transformations are hyperbolic transformations. Hence by using the Poincar\'e's theorem, we get that the group $\{\gamma_{1},\gamma_1',\gamma_{2},\gamma_{2}',...,\gamma_{n-1},\gamma_{n-1}',\gamma_{n},\gamma_{n}' \}$  generates a discrete subgroup of $PSL^{*}(2,\mathbb R)$. Since there is an overlapping issue on each generator between $\gamma_i$ and $\gamma'_i$, the group $<\gamma_{1},\gamma_{2},...,\gamma_{n-1},\gamma_{n}>$ generates a discrete subgroup of $PSL^*(2, \mathbb{R})$. Let, $\Gamma^3$ $=$ \{$<\gamma_{1},\gamma_{2},...,\gamma_{n-1},\gamma_{n}>$ $\cap$ $PSL(2,\mathbb R)\}$, be the subgroup of hyperbolic M\"{o}bius transformations in the upper half plane contained in the group $<\gamma_{1},\gamma_{2},...,\gamma_{n-1},\gamma_{n}>$. Then, $\Gamma^3$ is a Fuchsian group. This group $\Gamma^3$ is called the rank $n$, $n$ $\in$ $\mathbb{N} - \{1\}$, generalized Fuchsian Schottky group. Hence, the proof follows. 
                                 \end{proof}
                                
                                Therefore, by Proposition $4.1$, one can construct a rank $n$, $n$ $\in$ $\mathbb{N} - \{1\}$, generalized Fuchsian Schottky group by allowing any isometry of $\mathbb{H}$ to be a side-pairing transformation.\\

                                  The rest of this section is devoted to the study of the limit set for the Fuchsian Schottky group of rank $n$, $n \in \mathbb{N}$. Our constructed Fuchsian Schottky group of arbitrary finite rank $n$, $n \in \mathbb{N} - \{1\}$, is a particular example of a free Fuchsian group whose limit set is a Cantor set with sufficiently small limit set. The freeness of $\Gamma_{S_n}$ is defined by the well-known Ping-pong lemma. Though it is known that, a Kleinian group with a limit set of a Cantor set of Hausdorff dimension $< 1$ is always a free group. On the other hand, there exists a non-free purely hyperbolic Kleinian group with the limit set a Cantor set of Hausdorff dimension $< 1 + \epsilon$ (for $\epsilon > 0$) (see, \cite{Pankka}). Further, the limit points of this Fuchsian Schottky group $\Gamma_{S_n}$ lie within a line or a circle. Also, for the Fuchsian Schottky group, $\Gamma_{S_n}$, $n \in \mathbb{N}$, all the limit points are conical, i.e., each limit point of $\Gamma_{S_n}$ is a point of approximation of this group. So, the non-wandering set of the geodesic flow on the unitary tangent bundle of $\mathbb{H}/\Gamma_{S_n}$ is compact and none of the elements of this set diverges concerning the geodesic flow. Also, this set contains geodesic trajectories that are neither dense nor periodic. This characterizes the conical limit points in the limit set of $\Gamma_{S_n}$,  $n \in \mathbb{N}$. Hence, the group $\Gamma_{S_n}$ has no cusped limit point. Our constructed Schottky group satisfies the Beardon-Maskit condition in trivial, since 
                                  \{$\Lambda(\Gamma_{S_n})$ $-$  $\Lambda_{c}(\Gamma_{S_n})$\} $=$ $\phi$, $n$ $\in$ $\mathbb{N}$. Also, the measure, m\{$\Lambda(\Gamma_{S_n})$\} $=$ $0$ $=$  m\{$\Lambda_{c}(\Gamma_{S_n})$\}, which is the Ahlfors conjecture (general version one can say) for the group with the second kind. For $\Gamma_{S_1}$, the limit set consists of two points only. But for $\Gamma_{S_n}$, $n \in \mathbb{N} - \{1\}$, the limit set is a fractal with Hausdorff dimension $0< \delta = dim_H(\Lambda_{\Gamma_{S_n}}) < 1$. The ordinary set $\Omega(\Gamma_{S_n})$ has a single $\Gamma_{S_n}$-invariant component that is not simply connected. Since, for rank $n$ ($n$ $\in \mathbb{N})$, m\{$\Lambda_{c}(\Gamma_{S_n})$\} $\ne$ $1$, so the group $\Gamma_{S_n}$, $n$ $\in$ $\mathbb{N}$, is of convergence type and the Green functions exist on $\mathbb{H}$/$\Gamma_{S_n}$ concerning the hyperbolic Laplacian. The exponent of convergence of this group (denotes $\delta_{\Gamma_{S_n}}$) is bounded by $0$ and $\frac{1}{2}$, i.e., $0 < \delta_{\Gamma_{S_n}} \le \frac{1}{2}$. Since $\Gamma_{S_n}$ is geometrically finite, $\delta_{\Gamma_{S_n}}$ is the Hausdorff dimension of the limit set $\Lambda_{\Gamma_{S_n}}$.

                 \section{\textbf{ASSOCIATED SURFACES WITH THE CONFORMAL BOUNDARY AT INFINITY CORRESPONDING TO THE GROUP $\Gamma_{S_n}$ FROM THE POINT OF VIEW OF THE EULER CHARACTERISTIC, IN THE HYPERBOLIC SURFACE}}
                  
                  In this section, we have discussed the associated surfaces corresponding to rank $n$, ($n$ $\in$ $\mathbb{N}$) Fuchsian Schottky group. For this purpose, first, we have obtained the Nielsen region for the group  $\Gamma_{S_n}$, $n \in \mathbb{N} - \{1\}$ and characterized this region in $\overline{\mathbb{H}}$. After that, we discussed the equivalency of the convex core and the compact core for the group $\Gamma_{S_n}$. Further, we have presented surface decomposition and gluing for $\Gamma_{S_n}$ to create the required hyperbolic surface with the conformal boundary at infinity from the point of view of the Euler characteristic, $\chi$, in the hyperbolic surface.
                                   
                                  Due to the absence of elliptic elements in $\Gamma_{S_n}$, the quotient hyperbolic surface is smooth. If $\Gamma_{S_n}$ had elliptic elements then the quotient surface would be an Orbifold, with singularities (conical) corresponding to the elliptic fixed points. Now, given an isometry group, (let, $\Gamma_{S_n}$) of $\mathbb H $, a natural way to obtain a hyperbolic surface is as a quotient $\mathbb H$/$\Gamma_{S_n}$. Points in the quotient correspond to orbits of   $\Gamma_{S_n}$ and there is a natural projection $\Pi$: $\mathbb H $ $\to$ $\mathbb H$/$\Gamma_{S_n}$  given by $\Pi(z)$ =  $\Gamma_{S_n}(z)$. Since the action of $\Gamma_{S_n}$ is properly discontinuous on the hyperbolic plane, $\mathbb H $, the quotient is well-defined. Hence, the orbits are locally finite.\\

                                      5.1. \textbf{Surface decomposition for the group $\Gamma_{S_n}$, $n$ $\in \mathbb{N} - \{1\}$.}
                                      
                                       Here, $\Gamma_{S_n}$ is a non-elementary Fuchsian group of the second kind with $n$ number of generators. So, the corresponding hyperbolic surface ($X =$ $\mathbb H$/$\Gamma_{S_n}$) is non-elementary and geometrically finite. Also,  $\Lambda(\Gamma_{S_n})$ is nowhere dense in the circle at infinity and hence \{$\partial\mathbb H $ $-$  $\Lambda(\Gamma_{S_n})$\} is a countable union of open intervals $I_i$, where $i=1, 2, ... ,(2n-1), (2n), ...$ . Then the group $\Gamma_{S_n}$  acts properly discontinuously on \{$\mathbb{RP}^1$ $-$ $\Lambda(\Gamma_{S_n})$\}. Since $\Gamma_{S_n}$ is geometrically finite containing only hyperbolic elements as peripheral elements, the conformal boundary at infinity for $\Gamma_{S_n}$, i.e., \{$\mathbb{RP}^1$ $-$ (Cantor set with $dim_H \le \frac{1}{2}$ as a limit set of $\Gamma_{S_n}$)\}/ $\Gamma_{S_n}$ consists of a finite number of simple closed curves. The structure of the ends of $\mathbb{H}/\Gamma_{S_n}$ reflects the $1$-dimensional measure of $\Lambda(\Gamma_{S_n})$ with $\Lambda(\Gamma_{S_n})$ $\ne$ $\mathbb{RP}^1$ which gives the length of $\Lambda(\Gamma_{S_n})$ is zero. Now, let us consider that $\alpha_i$ is a geodesic whose endpoints are the endpoints of $I_i$ and $H_i$ is the half-plane bounded by $\alpha_i$ and $I_i$.
                                      Then, the Nielsen region (also called the convex hull of the limit set $\Lambda(\Gamma_{S_n})$ in $\mathbb{H}$) of rank $n$, $n$ $\in$ $\mathbb{N} - \{1\}$, Fuchsian Schottky group $\Gamma_{S_n}$ is defined by the subsequent set, 
                                      \begin{equation}
                                      N(\Gamma_{S_n}) = \{\mathbb H - (H_1 \cup H_2 \cup ... \cup H_{2n-1} \cup H_{(2n-1)+1} \cup ...)\}.
                                      \end{equation} The collection $\{(\cup H_i)$, $i = 1, 2, ..., (2n-1), (2n), ...\}$ is invariant under $\Gamma_{S_n}$ by the invariance of $\Lambda(\Gamma_{S_n})$. Since the cardinality of $\Lambda(\Gamma_{S_n})$ is greater than $1$, the convex core for $\Gamma_{S_n}$, i.e., $N(\Gamma_{S_n})$/$\Gamma_{S_n}$ is nonempty. Each semi-circular line, i.e., $D_i$ marks the boundary of the fundamental domain for the action of $\Gamma_{S_n}$. Since $\Gamma_{S_n}$ contains no parabolic element, there is zero possibility to create the truncated Nielsen region. If there is a parabolic element, say $c$ in $\Gamma_{S_n}$ then we will define the truncated Nielsen region as, ${N_1}(\Gamma_{S_n})$ $=$ $N(\Gamma_{S_n})$ $-$ $O_c$, where $O_c$ is the interior of horocycle $O'_c$ which will be created at the point $c$. So, the compact core (say, $K$) of $X$ is the same as the convex core in $X$, i.e., $C(X)$ = $K$. Hence, the compact core $K$ of the hyperbolic surface corresponding to $\Gamma_{S_n}$ is compact and $X-K$  is a finite disjoint union of funnels (only). Due to the absence of parabolic elements in $\Gamma_{S_n}$, the convex core boundary for $\Gamma_{S_n}$ doesn't contain any cusp end. Also, the convex core is compact, i.e., $\Gamma_{S_n}$ is a convex cocompact group. That means $\mathbb H$/$\Gamma_{S_n}$ is a conformally compact hyperbolic surface. Now, the Dirichlet fundamental domain, $D(\Gamma_{S_n})$ of $\Gamma_{S_n}$ has finitely many sides, so $D(\Gamma_{S_n})$ meets $\Omega(\Gamma_{S_n})$ in a finite collection of disjoint arcs \{$\alpha_i$\} which lie in the half-planes \{$H_1$, $H_2$, ..., $H_{2n-1}$, $H_{(2n-1)+1}$, ...\}. Therefore, the quotient $\overline{D(\Gamma_{S_n})\cap (\cup H_i)}/\Gamma_{S_n}$ (closure is taking in the Riemann-sphere topology) gives a finite collection of ends in the hyperbolic surface with an extra boundary circle at infinity. The collection of extra funnel end components corresponded with the conformal boundary of $\Gamma_{S_n}$. Therefore, $X$ has the decomposition as the compact core with only a finite number of funnels attached to it. So, we have performed the decomposition for purely hyperbolic generalized Schottky group as 
                                      \begin{equation}
                                      X = K \cup F'
                                      \end{equation} together with the funnels grouped  as the disjoint union, $\{F'_{1}$ $\cup$ $F'_2$ $\cup$ ... $\cup$ $ F'_k \}$ $=$ $F'$, where, $k$ is the number of funnel ends in the hyperbolic surface $X$ corresponding to the Fuchsian Schottky group of some finite rank. 
                                      
                                      Now, taking the notion from the preceding discussion, in the following, we have provided a characterization of the Nielsen region, $N(\Gamma_{S_n})$, for the group $\Gamma_{S_n}$, $n \in \mathbb{N} - \{1\}$, in $\overline{\mathbb{H}}$: 
                                                        
                                                        \textbf{(i)} Outside the Fuchsian Schottky curves on the hyperbolic plane, $\mathbb{H}$, $N(\Gamma_{S_n})$ of $\Gamma_{S_n}$, is bounded by the Euclidean half-planes (i.e., $H_1$, $H_2$, ..., $H_{2n}$, ...) centered on $\mathbb{R}$. Actually, $N(\Gamma_{S_n})$ is bounded by the axes of the hyperbolic elements in $\Gamma_{S_n}$ which project to the closed geodesics  and cut off the funnel ends on $\mathbb{H}/\Gamma_{S_n}$. 
                                                        
                                                        \textbf{(ii)} Inside the Fuchsian Schottky curves, the Nielsen region $N(\Gamma_{S_n})$ of $\Gamma_{S_n}
                                                        $ is unbounded on the hyperbolic plane, $\mathbb{H}$ at $\mathbb{RP}^1$.\\
                                                        
                                                        5.2. \textbf{The conformal boundary at infinity corresponds to the group $\Gamma_{S_n}, n \in \mathbb{N}$, in the point of view of the topological invariant Euler characteristic, in the hyperbolic surface. } 
                                                                          
                                                                          The region $F_{S_n}$ $=$ \{$\mathbb H $ $-$ $( D_1 \cup D_2 \cup ... \cup D_{2n-1} \cup D_{2n})$\}  is the fundamental domain for the action of $\Gamma_{S_n}$. It is well known that $\mathbb H$/$\Gamma_{S_n}$ and ${\overline{F}_{S_n}}$/$\Gamma_{S_n}$ are topologically equivalent  if and only if $F_{S_n}$ is locally finite (see, \cite{Beardon}) (${\overline{F}_{S_n}}$ is the closure of $F_{S_n}$). $F_{S_n}$ is a geodesically convex fundamental domain, so it is also a Dirichlet domain. Again we know that the Dirichlet domain is always locally finite. Note that, the fundamental region $F_{S_n}$ has Euler characteristic, $\chi$ $=$ $1$. On the other hand, for any hyperbolic surfaces, the Euler characteristic is 
                                                                          \begin{equation}
                                                                          \chi (X) = 2- 2g- n_{F'}- n_{C}
                                                                          \end{equation} 
                                                                          where $g$ is the number of the genus, $n_{F'}$ denotes the number of funnel ends and $n_C$  indicates the number of cusp ends. Now, for the group $\Gamma_{S_n}$ after gluing $2n$ edges together to form the required surface $X$, we obtain the Euler characteristic 
                                                                          \begin{equation}
                                                                          \chi = (1-n)
                                                                          \end{equation}
                                                                          \\Therefore, combining equations $(5.3)$ and $(5.4)$ we get the following $\textbf{(}$keeping equation $(5.2)$ in mind$\textbf{)}$:

                                                                           For rank, $n$=$1$:
                                                                           \begin{equation}
                                                                           n_{F'} + 2g = 2
                                                                           \end{equation}
                                                                           So, only one case arises; $g=0, n_{F'}=2$.
                                                                           
                                                                          For rank, $n$=$2$: 
                                                                          \begin{equation}
                                                                          n_{F'} + 2g = 3 
                                                                          \end{equation}
                                                                          Therefore, $2$ cases arise; (i) $g=1, n_{F'} = 1$, and (ii) $g=0, n_{F'} =3$.
                                                                          
                                                                          For rank, $n$=$3$: 
                                                                          \begin{equation}
                                                                          n_{F'} + 2g = 4
                                                                          \end{equation}
                                                                          That means, here also $2$ cases arise; (i) $g=1, n_{F'} =2$, and (ii) $g=0, n_{F'}=4$.
                                                                          
                                                                          For rank, $n$=$4$:
                                                                          \begin{equation}
                                                                          n_{F'} + 2g = 5
                                                                          \end{equation}
                                                                          Hence, $3$ cases arise; (i) $g=1, n_{F'}=3$, (ii) $g=2, n_{F'}=1$, and (iii) $g=0, n_{F'} =5$.
                                                                          
                                                                          For rank, $n$=$5$:
                                                                          \begin{equation}
                                                                           n_{F'} + 2g = 6
                                                                          \end{equation}
                                                                          That means, here also $3$ cases arise; (i) $g=0, n_{F'} = 6$, (ii) $g=1, n_{F'} =4$, and (iii) $g=2, n_{F'} =2$.
                                                                          
                                                                          For rank, $n$=$6$:
                                                                          \begin{equation}
                                                                          n_{F'} + 2g = 7
                                                                          \end{equation}
                                                                          So, $4$ cases arise; (i) $g=0, n_{F'}=7$, (ii) $g=1, n_{F'}=5$, (iii) $g=2, n_{F'}=3$, and (iv) $g=3, n_{F'} =1$.
                                                                          
                                                                          For rank, $n$=$7$:
                                                                          \begin{equation}
                                                                          n_{F'} + 2g = 8
                                                                          \end{equation}
                                                                          i.e., here also $4$ cases arise; (i) $g=0, n_{F'}=8$, (ii) $g=1, n_{F'}=6$, (iii) $g=2, n_{F'} =4$, and (iv) $g=3, n_{F'} =2$.
                                                                          
                                                                          For rank, $n$=$8$:
                                                                          \begin{equation}
                                                                          n_{F'} + 2g = 9
                                                                          \end{equation}
                                                                          Therefore, $5$ cases arise; (i) $g=0, n_{F'}=9$, (ii) $g=1, n_{F'}=7$, (iii) $g=2, n_{F'}=5$, (iv)$g=3, n_{F'}=3$, and (v) $g=4, n_{F'}=1$.
                                                                          
                                                                          For rank, $n$=$9$:
                                                                          \begin{equation}
                                                                                            n_{F'} + 2g = 10
                                                                                            \end{equation}
                                                                          Hence, $5$ cases again arise: (i) $g=0, n_{F'}=10$, (ii) $g=1, n_{F'}=8$, (iii) $g=2, n_{F'}=6$, (iv) $g=3, n_{F'}=4$, and (v) $g=4, n_{F'}=2$.
                                                                          
                                                                           For rank, $n$=$10$:
                                                                                            \begin{equation}
                                                                                                              n_{F'} + 2g = 11
                                                                                                              \end{equation}
                                                                                            So, $6$ cases arise: (i) $g=0, n_{F'}=11$, (ii) $g=1, n_{F'}=9$, (iii) $g=2, n_{F'}=7$, (iv) $g=3, n_{F'}=5$, (v) $g=4, n_{F'}=3$, (vi) $g=5, n_{F'}=1$;\\
                                                                                            
                                                                                             and so on. \\
                                                                                            
                                                                            \textbf{Observations:} 
                                                                            
                                                                            For $\Gamma_{S_2}$: Equation $(5.6)$(i) gives a hyperbolic surface containing $1$ funnel end with $1$ genus. Equation $(5.6)$(ii) produces a hyperbolic surface containing $3$ funnel ends with no genus. So, for rank $2$ the maximum number of funnel ends is $3$. Notice that, one can attach two funnel ends with a suitable twist parameter and produce a genus $0$ surface. So, equation $(5.6)$(ii) surface can be converted to equation $(5.6)$(i) surface.
                                                                          
                                                        	                  For $\Gamma_{S_3}$: Equation $(5.7)$(i) produces a surface containing $2$ funnel ends with $1$ genus. Equation $(5.7)$(ii) gives a surface containing $4$ funnel ends with no genus. So, here also attaching two funnel ends with a suitable twist parameter we get a surface containing $1$ genus and $2$ funnel ends. Note that, keeping in mind equation $(5.2)$, we can not go further to get $2$ genus surfaces with no funnel end. Observe that, for the rank $3$ Fuchsian Schottky group, the maximum number of funnel ends can be $4$. 
                                                                          
                                                                          For $\Gamma_{S_4}$: Equation $(5.8)$(i) delivers a surface containing $3$ funnel ends with $1$ genus. Now using a suitable twist parameter we get equation $(5.8)$(ii), i.e., the surface containing $2$ genera with $1$ funnel end. Further, Equation $(5.8)$(iii) produces a surface containing $5$ funnel ends (which is the maximum for $\Gamma_{S_4}$) with no genus. Observe that, this surface can be converted to the equation  $(5.8)$(ii) surface by using two suitable twist parameters. 
                                                                          
                                                                          For $\Gamma_{S_5}$: Equation $(5.9)$(i) produces a surface containing no genus with $6$ funnel ends (which is the maximum for $\Gamma_{S_5}$). Notice that, surfaces coming from $(5.9)$(i) and $(5.9)$(ii) can be reduced in the equation $(5.9)$(iii) surface by using $2$ and $1$ suitable twist parameters respectively. 
                                                                          
                                                                          For $\Gamma_{S_6}$: Equation $(5.10)$(i) gives a surface containing no genus with $7$ funnel ends (which is the maximum for $\Gamma_{S_6}$). Observe that, surfaces arising from $(5.10)$(i), $(5.10)$(ii), and $(5.10)$(iii) can be reduced in the equation $(5.10)$(iv) surface by using $3$, $2$, and $1$ suitable twist parameters respectively.  
                                                                          
                                                                          For $\Gamma_{S_7}$: Equation $(5.11)$(i) delivers a surface containing no genus with $8$ funnel ends (which is the maximum for $\Gamma_{S_7}$). Note that, surfaces originating from $(5.11)$(i), $(5.11)$(ii), and $(5.11)$(iii) can be converted in the equation $(5.11)$(iv) surface by using $3$, $2$, and $1$ suitable twist parameters respectively.
                                                                          
                                                                          For $\Gamma_{S_8}$: Equation $(5.12)$(i) gives a surface containing no genus with $9$ funnel ends (which is the maximum for $\Gamma_{S_8}$). Observe that, surfaces coming from $(5.12)$(i), $(5.12)$(ii), $(5.12)$(iii), and $(5.12)$(iv) can be consolidated in the equation $(5.12)$(v) surface by using $4$, $3$, $2$, and $1$ suitable twist parameters respectively.
                                                                          
                                                                          For $\Gamma_{S_9}$: Equation $(5.13)$(i) produces a surface containing no genus with $10$ funnel ends (which is the maximum for $\Gamma_{S_9}$). Notice that, surfaces arising from $(5.13)$(i), $(5.13)$(ii), $(5.13)$(iii), and $(5.13)$(iv) can be reduced in the equation $(5.13)$(v) surface by using $4$, $3$, $2$, and $1$ suitable twist parameters respectively.
                                                                          
                                                                          For $\Gamma_{S_{10}}$: Equation $(5.14)$(i) delivers a surface containing no genus with $11$ funnel ends (which is the maximum for $\Gamma_{S_{10}}$). Note that, surfaces coming from $(5.14)$(i), $(5.14)$(ii), $(5.14)$(iii), $(5.14)$(iv), and $(5.14)$(v) can be converted in the equation $(5.14)$(vi) surface by using $5$, $4$, $3$, $2$, and $1$ suitable twist parameters respectively; \\
                                                                          
                                                                          and so on.\\

                                                                          \textbf{Conclusion of the overhead observations:} 
                                                                          
                                                                          (i) For the rank $n$ Fuchsian Schottky group, the associated hyperbolic surface contains the maximum number of funnel ends when the surface doesn't contain a genus. In this case, the number of funnels end is equal to $(n+1)$.
                                                                          
                                                                          (ii) If we apply the attaching process with the funnel ends by using suitable twist parameters then we get, $n=2g$, when  $ n_f=1$ and $n= 2g+ 1 $, when $ n_f= 2$ ; where $n \in \mathbb{N} - \{1\}$ and $g \in \mathbb{N}$.
                                                                          
                                                                           So, the preceding conclusion leads to the following two theorems:
                                                                                             
                                                                                             \begin{thm} \label{t5.1}
                      For the Euler characteristic, $\chi$, in the hyperbolic surface, the quotient surface corresponding to the rank $n$ Fuchsian Schottky group 
                                                                                                  
                                                                                                  (i) is of $\frac{(n+2)}{2}$ (when $n$ is even) or $\frac{(n+1)}{2}$ (when $n$ is odd) types of hyperbolic surfaces.
                                                                                                  
                                                                                                  (ii) contains the maximum number of funnels end (n+1) when the surface has no genus.
                                                                                                  
                                                                                                  (iii) contains the minimum number of funnel ends $1$ or $2$ when the surface has genus $\frac{n}{2}$ (for $n$ is even) or $\frac{(n-1)}{2}$ (for $n$ is odd) respectively.                                                                       
                                                                                             \end{thm}

                                                                                            \begin{thm} \label{t5.2}
                                                                                            For the Euler characteristic, $\chi$, in the hyperbolic surface, the associated hyperbolic surface (applying attaching process with suitable twist parameters) corresponding to the Fuchsian Schottky group, with 
                                                                                                                                                                       
                                                                                                                                                                       (i) rank $1$, contains two funnel ends with no genus, i.e., it is a hyperbolic cylinder.

                                                                                                                                                                        (ii) even rank (for some even number $m, m$ $\in$ $2 \mathbb{N}$), contains $\frac{m}{2}$ genera with one funnel end in the conformal boundary at infinity; i.e., it is the finite Loch Ness monster, having arbitrarily finitely many handles and only one way to go to infinity.
                                                                                                                                                                       
                                                                                                                                                                       (iii) odd rank \{for some odd number $n$, $n$ $\in$ $(2r-1)\mathbb{N}$, ($r \in \mathbb{N}$, but $r$ $\ne 1)$\}, contains $\frac{(n -1)}{2}$ genera with two funnels ending in the conformal boundary at infinity; i.e., it is the finite Jacob's ladder, having two ways to go to infinity and arbitrarily finitely many handles in each.
                                                                                            \end{thm} 
                                                                                            
                                                                                            \section{\textbf{THE PAIR OF PANTS SITUATION FOR THE GROUP $\Gamma_{S_n}$, $n \in \mathbb{N} - \{1\}$}}
                                                                                           
                                                                                            \textbf{Existence of a non-tight pair of pants for $\Gamma_{S_n}$, $n$ $\in$ $\mathbb{N}$ with $n$ $\ne$ $1$.} 
                                                                                                                                
                                                                                                                                Assume that, $\gamma$ $\in$ $\Gamma_{S_n}$ $-$ \{Id\} and $a_0$ $\in$ $\mathbb{H}$ such that $a_0$ is not fixed by $\gamma$. So, the perpendicular bisector of the hyperbolic segment [$a_0,\gamma(a_0)$] defined by $P_{a_0}(\gamma)$ $=$ $\{a \in \mathbb{H} : d(a,a_0) = d(a,\gamma(a_0))\} $ is the geodesic that separates $\mathbb{H}$ into two connected components, orthogonal to the segment [$a_0,\gamma(a_0)$], passing through its middle point. Let, $C_{a_0}$ be the closed half-plane in $\mathbb{H}$ bounded by $P_{a_0}(\gamma)$ containing $\gamma(a_0)$. So, $\gamma(P_{a_0}(\gamma^{-1}))$ $=$ $P_{a_0}(\gamma)$ and $\gamma(C^{\circ}_{a_0}(\gamma ^{-1}))$ $=$ $\mathbb{H}$ $-$ $C_{a_0}(\gamma)$ hold, where $C^{\circ}(\gamma)$ denotes the interior of $C(\gamma)$. Since, $\Gamma_{S_n}$ is purely hyperbolic, so for an element $\gamma$ in $\Gamma_{S_n}$,  $\overline{P}_{a_0}(\gamma)$ $\cap$ $\overline{P}_{a_0}(\gamma^{-1})$ $=$ $\phi$ holds in $\overline{\mathbb{H}}$. Now, for $\gamma$ $\in$ $\Gamma_{S_n}$,  $C(\gamma)$ $\cap$ $C(\gamma ^{-1})$ $=$ $\phi$ and  $\gamma(C(\gamma ^ {-1}))$ $=$ $\mathbb{H}$ $-$ $C^{\circ}(\gamma)$. Hence, only two types of cases can arise here for the group $\Gamma_{S_n}$. For $\gamma_1$,$\gamma_2$, ...,$\gamma_n$ $\in$ $\Gamma_{S_n}$, the sets \{$C(\gamma_1)$ $\cup$ $C(\gamma_1^{-1})$\},  \{$C(\gamma_2)$ $\cup$ $C(\gamma_2^{-1})$\}, ..., \{$C(\gamma_n)$ $\cup$ $C(\gamma_n^{-1})$\} are disjoint.
                                                                                                                                 So, we have two cases as in the following.
                                                                                                                                 
                                                                                                                                 \textbf{(i)} \{$C(\gamma_i)$ $\cap$ $C(\gamma_i ^{-1})$\} $=$ $\phi$ for $\gamma_i$ $\in$ $\Gamma_{S_n}$ $(i= 1,2, ... ,n)$ and they are adjacent to each of them in $\mathbb{H}$.
                                                                                                                                 
                                                                                                                                 \textbf{(ii)} \{$C(\gamma_j)$ $\cap$ $C(\gamma_j ^ {-1})$\} $=$ $\phi$ for $\gamma_j$ $\in$ $\Gamma_{S_n}$ $(j=1,2, ... ,n)$ and they are not adjacent, rather opposite to each of them in $\mathbb{H}$.
                                                                                                                                  \\ Clearly, from our construction of  rank $n$, $n$ $\in$ $\mathbb{N} - \{1\}$, Fuchsian Schottky group in Proposition $3.1$, we claim that only case (i) can occur for the group $\Gamma_{S_n}$, $n$ $\in$ $\mathbb{N} - \{1\}$ because we have chosen the semi-circles with diameters lying on the real axis and with their reflection in the upper imaginary axis. Now, we define the group $\Gamma_{S_n}$ is called crossed when we order the sets $C(\gamma_i^{\pm 1})$, $(i=1,2,...,n)$ along the boundary of $\mathbb{H}$. Observe  that, $C(\gamma_i)$ and $C(\gamma_i ^{-1})$ are not adjacent. So, any rank $n$, $n$ $\in$ $\mathbb{N}- \{1$\}, Fuchsian Schottky group  is not crossed. Hence, the non-tight pair of pants arises from the group  $\Gamma_{S_n}$, which we have described below in brief.\\

                         \textbf{Attaching $Y$-pieces for $\Gamma_{S_2}$.} 
                                                                                                                                  
                                                                                                                                  The Fuchsian Schottky group of rank $2$ produces pairs of pants with no cusp holes where $\gamma_1$ and $\gamma_2$ are hyperbolic with $\gamma_1\gamma_2$ hyperbolic. For hyperbolic elements $\gamma_1$, $\gamma_2$, and $\gamma_1$$\gamma_2$ to match with the conformal holes via their fixed points, we have the axes $A(\gamma_1)$, $A(\gamma_2)$, and $A(\gamma_1\gamma_2)$ to connect its fixed points by the geodesic to bound the hole. Now, taking the notion of the previous section, we obtain for $\Gamma_{S_2}$ that the quotient $\mathbb{H}/\Gamma_{S_2}$ gives a hyperbolic handlebody of genera $2$. So, by pants decomposition, for a rank $2$ Fuchsian Schottky group, it provides a closed orientable hyperbolic surface which can be subdivided by $3$ closed, simple, disjoint geodesics into the union of $2$ pairs of pants with the geodesics of their boundary circles. Let, $Y^1$ and $Y^2$ be two $Y$-pieces with boundary geodesics $\gamma_1^1$, $\gamma_1^2$, $\gamma_1^1\gamma_2^1$ ; $\gamma_2^1$, $\gamma_2^2$,  and $\gamma_1^2\gamma_2^2$  parametrized on $\mathbb{R}$$/[t\rightarrow t+1]$, $t$ $\in [0,1]$. Note that, here $\{\gamma_1^1,\gamma_2^1,\gamma_1^1\gamma_2^1 \}$ and $\{\gamma_1^2, \gamma_2^2, \gamma_1^2\gamma_2^2\}$ all are hyperbolic. So, the lengths of these boundary geodesics are always positive. Let us suppose that $\gamma_1^1\gamma_2^1$ $=$ $\gamma_3^1$  and $\gamma_1^2\gamma_2^2$ $=$ $\gamma_3^2$. Now, $l(\gamma_i^j)$ = $d(z_k,\gamma_i^j(z))$, where $z_k$ is the axis of $\gamma_i^j$ $(i=1,2,3$ ; $j=1,2$ ; $ k=1,2,3)$. Further, assume that,  $l(\gamma_1^1)$ $=$ $l(\gamma_2^2)$. Then for any $\beta \in \mathbb{R}$ (call the twist parameter), we have built an $X$-piece via the identification $\gamma_1^1(t)$ $=$ $\gamma_2^2(\beta - t)$ $=$ $\gamma^\beta(t)$, where $t\in \mathbb{R}$$/[t\rightarrow t+1]$, $t\in [0,1]$. So, by pasting the two $Y$-pieces we obtain $X^\beta$ = $Y^1$ + $Y^2$ with the preceding twist parameter condition. Hence, for $\Gamma_{S_2}$, the corresponding hyperbolic surface is a conformally compact Riemann surface of signature $(0,4)$, i.e., it is a one $X$-piece whose length of all the boundary geodesics \{$\gamma_2^1, \gamma_3^1, \gamma_1^2, \gamma_3^2$\} are positive.\\ 
                                                                                                                                  
                                                                                                                                  \textbf{Attaching $Y$-pieces for $\Gamma_{S_3}$.} 
                                                                                                                                  
                                                                                                                                  For $\Gamma_{S_3}$, the quotient $\mathbb{H}/\Gamma_{S_3}$ contains $6$ boundaries, and for hyperbolic elements  $\gamma_1$, $\gamma_2$, $\gamma_3$, $\gamma_1\gamma_2$, $\gamma_2\gamma_3$, $\gamma_1\gamma_3$ that match with the conformal holes via their fixed points. Hence, for $\Gamma_{S_3}$ the closed orientable hyperbolic surface can be subdivided by $6$ simple, closed, disjoint geodesics into the union of $4$ pairs of pants with the geodesics of their boundary circles. Let, $Y^1$, $Y^2$, $Y^3$ and $Y^4$ be four $Y$-pieces with boundary geodesics $\gamma_1^1, \gamma_2^1, \gamma_1^1\gamma_2^1$ ; $\gamma_2^2, \gamma_3^2, \gamma_2^2\gamma_3^2$ ; $\gamma_1^3, \gamma_3^3, \gamma_1^3\gamma_3^3$ ; $\gamma_1^4\gamma_2^4, \gamma_2^4\gamma_3^4, \gamma_1^4\gamma_3^4$ respectively parametrized on $\mathbb{R}/[t\rightarrow t+1]$, $t \in [0,1]$. Let, $\gamma_1^1\gamma_2^1$ $=$ $\gamma_3^1$, $\gamma_2^2\gamma_3^2$ $=$ $\gamma_5^2$, $\gamma_1^3\gamma_3^3$ $=$ $\gamma_4^3$, $\gamma_1^4\gamma_2^4$ $=$ $\gamma_3^4$, $\gamma_2^4\gamma_3^4$ $=$ $\gamma_5^4$, $\gamma_1^4\gamma_3^4$ $=$ $\gamma_4^4$. Notice that, $l(\gamma_i^j)$ $=$ $d(z_k,\gamma_i^j(z))$, where $z_k$ is the axis of $\gamma_i^j$ $(i=1,2,3,4,5$ ; $j=1,2,3,4$ ; $ k=1,2,3,4,5)$ and lengths of all of these elements are positive. Let, $l(\gamma_1^1)$ $=$ $l(\gamma_2^2)$ and $l(\gamma_3^3)$ $=$ $l(\gamma_4^4)$. Then for twist parameters $\beta^1$, $\beta^2$ $\in$ $\mathbb{R}$, we prepare two $X$-pieces via the  identifications $\gamma_1^1(t)$ $=$ $\gamma_2^2(\beta^1 - t)$ $=$ $\gamma^{\beta^1} (t)$ and $\gamma_3^3(t)$ $=$ $\gamma_4^4(\beta^2 - t)$ $=$ $\gamma^{\beta^2} (t)$, where $t\in \mathbb{R}/[t\rightarrow t+1],$ $t\in [0,1].$ Now, by pasting these four $Y$-pieces we get, $X^{\beta^1}$ $+$ $X^{\beta^2}$ $=$ ($Y^1$ $+$ $Y^2$) $+$ ($Y^3$ $+$ $Y^4$) $=$ $X^1$ $+$ $X^2$,  with the overhead two twist parameter conditions. Now, suppose that $l(\gamma_3^2)$ $=$ $l(\gamma_1^3)$ and $l(\gamma_5^2)$ $=$ $l(\gamma_4^3)$. Then for twist parameters $\beta^3$, $\beta^4$  $\in$ $\mathbb{R}$, we have finally built a hyperbolic surface $X^3$ of signature $(1,4)$ via the two identifications $\gamma_3^2 (t)$ $=$ $\gamma_1^3 (\beta^3 - t)$ $=$ $\gamma^{\beta^3} (t)$ and $\gamma_5^2 (t)$ $=$ $\gamma_4^3 (\beta^4 - t)$ $=$ $\gamma^{\beta^4} (t)$ where $t\in \mathbb{R}/[t \rightarrow t+1]$, $t\in [0,1].$ Hence, for $\Gamma_{S_3}$ the corresponding hyperbolic surface is a conformally compact Riemann surface of signature $(1,4)$.\\
                                                                                                                                  
                                                                                                                                  \textbf{Attaching $Y$-pieces for $\Gamma_{S_4}.$ }
                                                                                                                                                                                        
                                                                                                                                  For $\Gamma_{S_4}$ the closed orientable hyperbolic surface can be subdivided by $9$ simple, closed, disjoint geodesics into the union of $6$ pairs of pants with the geodesics of their boundary circles. Let, $Y^1$, $Y^2$, $Y^3$, $Y^4$, $Y^5$ and $Y^6$ be six $Y$-pieces with boundary geodesics $\gamma_1^1, \gamma_2^1\gamma_3^1, \gamma_1^1\gamma_2^1\gamma_3^1$ ; $\gamma_2^2, \gamma_3^2\gamma_4^2, \gamma_2^2\gamma_3^2\gamma_4^2$ ; $\gamma_3^3, \gamma_1^3\gamma_4^3, \gamma_1^3\gamma_3^3\gamma_4^3$ ; $\gamma_4^4, \gamma_1^4\gamma_2^4, \gamma_1^4\gamma_2^4\gamma_4^4$ ; $\gamma_1^5\gamma_2^5, \gamma_3^5\gamma_4^5, \gamma_2^5\gamma_3^5$ ; and $\gamma_1^6\gamma_3^6, \gamma_2^6\gamma_4^6, \gamma_1^6\gamma_4^6$ respectively parametrized on $\mathbb{R}/[t \rightarrow t+1]$, $t$ $\in$ $[0,1]$. Let, $\gamma_2^1\gamma_3^1$ $=$ $\gamma_5^1$, $\gamma_1^1\gamma_2^1\gamma_3^1$ $=$ $\gamma_6^1$, $\gamma_3^2\gamma_4^2$ $=$ $\gamma_7^2$, $\gamma_2^2\gamma_3^2\gamma_4^2$ $=$ $\gamma_9^2$, $\gamma_1^3\gamma_4^3$ $=$ $\gamma_5^3$, $\gamma_1^3\gamma_3^3\gamma_4^3$ $=$ $\gamma_8^3$, $\gamma_1^4\gamma_2^4$ $=$ $\gamma_3^4$, $\gamma_1^4\gamma_2^4\gamma_4^4$ $=$ $\gamma_7^4$, $\gamma_1^5\gamma_2^5$ $=$ $\gamma_3^5$, $\gamma_3^5\gamma_4^5$ $=$ $\gamma_7^5$, $\gamma_2^5\gamma_3^5$ $=$ $\gamma_5^5$, $\gamma_1^6\gamma_3^6$ $=$ $\gamma_4^6$, $\gamma_2^6\gamma_4^6$ $=$ $\gamma_6^6$, $\gamma_1^6\gamma_4^6$ $=$ $\gamma_5^6$. Also, $l(\gamma_i^j)$ $=$ $d(z_k,\gamma_i^j(z))$, where $z_k$ is the axis of $\gamma_i^j$ $(i=1,2,3,4,5,6,7$ ; $j=1,2,3,4,5,6$ ; $ k=1,2,3,4,5,6,7)$. Assume that, $l(\gamma_1^1)$ $=$ $l(\gamma_2^2)$, $l(\gamma_3^3)$ $=$ $l(\gamma_4^4)$, $l(\gamma_5^5)$ $=$ $l(\gamma_6^6)$. Then for twist parameters $\beta^1$, $\beta^2$, $\beta^3$ $\in$ $\mathbb{R}$, we obtain three $X$-pieces via the two  identifications $\gamma_1^1 (t)$ $=$ $\gamma_2^2 (\beta^1 - t)$ $=$ $\gamma^{\beta^1}(t)$, $\gamma_3^3(t)$ $=$ $\gamma_4^4 (\beta^2 - t)$ $=$ $\gamma^{\beta^2} (t)$, and $\gamma_5^5$ $=$ $\gamma_6^6 (\beta^3 - t)$, where $t\in \mathbb{R}/[t\rightarrow t+1],$ $t\in [0,1].$ So, by pasting the six $Y$-pieces we get, $X^{\beta ^1}$ $+$ $X^{\beta^{2}}$ $+$ $X^{\beta^3}$ $=$ ($Y^1 + Y^2$) $+$ ($Y^3 + Y^4$) $+$ ($Y^5 + Y^6$) $=$ $X^1$ $+$ $X^2$ $+$ $X^3$ (say), with the exceeding three twist parameter conditions. Now, for $X^1$ and $X^2$, let $l(\gamma_7^2)$ $=$ $l(\gamma_5^3)$ and $l(\gamma_9^2)$ $=$ $l(\gamma_8^3)$. Then for twist parameters $\beta^4$, $\beta^5$  $\in$ $\mathbb{R}$, we have built a hyperbolic surface $X^{12}$ (say, $X^4$) of signature $(1,4)$ via the two identifications $\gamma_7^2 (t)$ $=$ $\gamma_5^3(\beta^4 - t)$ $=$ $\gamma^{\beta^4} (t)$ and  $\gamma_9^2(t)$ $=$ $\gamma_8^3 (\beta^5 - t)$ $=$ $\gamma^{\beta^5}(t)$, where $t$ $\in$ $\mathbb{R}/[t\rightarrow t+1]$, $t$ $\in$ $[0,1]$. Again, for $X^4$ and $X^3$, suppose $l(\gamma_3^4)$ $=$ $l(\gamma_3^5)$ and $l(\gamma_7^4)$ $=$ $l(\gamma_7^5)$. Then for twist parameters $\beta^6$, $\beta^7$ $\in$ $\mathbb{R}$, we have finally reached the required hyperbolic surface $X^{123}$ (say, $X^5$) of signature $(2,4)$ via the two identifications $\gamma_3^4 (t)$ $=$ $\gamma_3^5 (\beta^6 - t)$ $=$ $\gamma^{\beta^6} (t)$ and $\gamma_7^4 (t)$ $=$ $\gamma_7^5 (\beta^7 - t)$ $=$ $\gamma^{\beta^7} (t)$, where $t \in \mathbb{R}/[t \rightarrow t+1]$, $t \in [0,1]$. Hence, for $\Gamma_{S_4}$, the corresponding hyperbolic surface is a conformally compact Riemann surface of signature $(2,4)$.
                                                                                                                                  
                                                                                                                                 In this way, for the rank $n$ Fuchsian Schottky group $\Gamma_{S_n}$, $n \in \mathbb{N} - \{1\}$, one can attach $(2n - 2)$ numbers of non-tight pairs of pants by using  $(3n-2)$ numbers of twist parameters where Bers' constant $<$ $(31n + 21)$ (for Bers' constant the readers may go through the book of Buser \cite{Buser}) to create the required associated hyperbolic surface $X^*$ (say), which is basically the conformally compact Riemann surface of signature $\textbf{(}(n-2), 4\textbf{)}$.\\
                                                                                                                                 
                                                                                                                                \textbf{Half-collars for $\Gamma_{S_n}$.} 
                                                                                                                                
                                                                                                                                For pair of pants $Y^k$, $k=1,2, ..., (2n-2)$ with boundary geodesics $\gamma^k_i$, $i=1,2,3$, the sets $ \mu[\gamma_i^k]$ $=$ $\{a \in Y^k, a \notin \gamma^k_i : 0 < \sinh(dist(a, \gamma_i^k)) \sinh \frac{1}{2}l(\gamma_i^k) \le 1\}$ are homeomorphic to \{$(0, 1]$ $\times$  $\mathbb{R}/[t\rightarrow t+1]$\} with the funnel metric. Note that, here dist\{$\mu[\gamma^k_i], \mu[\gamma^k_j]\} \nrightarrow 0$ for any  $i\ne j$. Consequently, the sets $\mu[\gamma^k_i]$ are always pairwise disjoint for $\Gamma_{S_n}$. So, the collar around $\gamma^k_i$ for the group $\Gamma_{S_n}$ is defined by the subsequent set, $$C(\gamma^k_i) = \{a\in X^* : 0 < dist(a,\gamma^k_i) \le w(\gamma^k_i)\},$$ where, $w(\gamma^k_i) =$ arcsinh$[\frac{1}{\sinh(\frac{1}{2}l(\gamma^k_i))}]$ is the width. Since there is no thin part of the hyperbolic surface $X^*$, the $width$ of the half-collars  is minimal for the group $\Gamma_{S_n}$. This characterizes any arbitrary finite rank Fuchsian Schottky group from the point of view of the Fuchsian surface group in the hyperbolic space.\\

                                                                                                                                \textbf{Fenchel-Nielsen coordinates for the group $\Gamma_{S_n}$.} 
                                                                                                                                
                                                                                                                                The deformation parameters for the hyperbolic surface $X^*$ corresponding to rank $n$ Fuchsian Schottky group obtained through our preceding study are $l_1, l_2, ..., l_{(3n-2) + 4}$ and $\beta_1, \beta_2, ..., \beta_{(3n-2)}$. Hence, the Fenchel-Nielsen coordinates of Teichm\"uller space for the group $\Gamma_{S_n}$ which is a subset of real $(6n-4)$ dimensional Euclidean space, $\mathbb{R}^{6n-4}$ are given by : $$ T_{X^*} = \{l_1, l_2, ..., l_{(3n+1)}, l_{(3n+2)} ; \beta_1, \beta_2, ..., \beta_{(3n-3)}, \beta_{(3n-2)} \in \mathbb{R}_+^{(3n+2)} \times \mathbb{R}^{(3n-2)}\}$$ 
                                                                                                                                  Therefore, from the above discussion we can state the following: 
                                                                                                                                  \begin{thm}
                                                                                                                                  \label{t6.1}
                                                                                                                                                                        $\textbf{(Non-tight pants decomposition for any finite rank Fuchsian}$\\
                                                                                                                                                                        $\textbf{ Schottky group)}$ \\ The compact core of the conformally compact Riemann surface $X^*$ corresponding to the Fuchsian Schottky group of rank $n$, $n \in \mathbb{N} - \{1\}$, can be decomposed into $2(n-1)$ numbers of non-tight pairs of pants $Y_i$ by using $(3n-2)$ numbers of twist parameters where Bers' constant is strictly less than $(31n+21)$, such that $X^* = \{Y_1 \cup Y_2 \cup ... \cup Y_{2(n-1)}\} \cup \{F'_1 \cup F'_2 \cup F'_3 \cup F'_4\}$ with $(n-2)$ genera, where $F'_j$ denotes the funnel ends, $j=1,2,3,4$, and $(2-2n)$ is the Euler characteristic, $\chi$ of $X^*$. Moreover, the Fenchel-Nielsen coordinates for  Teichm\"uller space ($\subset \mathbb{R}^{6n-4}$) corresponding to the group $\Gamma_{S_n}$ are $\{l_1, l_2, ..., l_{(3n+2)}$ ; $\beta_1, \beta_2, ..., \beta_{(3n-2)}\}$.
                                                                                                                                  \end{thm}
                 
                 \section {\textbf{CONCLUSION}}
                  Though any finite rank Fuchsian Schottky group is a discrete subgroup of $PSL(2, \mathbb{R})$, in this paper, we have introduced orientation-reversing isometries as side-pairing transformations that enrich our constructed group from the group's theoretical point of view (see, Proposition $4.1$). We noticed that the convex cocompact Fuchsian Schottky hyperbolic surface had a sufficiently small limit set with Hausdorff dimension $0 \le \delta_{\Gamma_{S_n}} \le \frac{1}{2}$. We also observed that, from the point of view of the topological invariant Euler characteristics in the hyperbolic surface, the Fuchsian Schottky group produced hyperbolic surfaces like finite Loch Ness monster and finite Jacob's ladder (see, Theorem $5.1$ and Theorem $5.2$). Note that, for a compact Riemann surface corresponding to a Fuchsian group of genus $n$ $(n>1)$ and a conformally compact Riemann surface corresponding to the Fuchsian group equipped with the Schottky structure of rank $n$ $(n>1)$ (which we have constructed in this paper), the dimensions of Teichm\"uller spaces are not equal (increased by $2$ in the sense of real dimensional Euclidean space for $\Gamma_{S_n}$, see, Theorem $6.1$).\\

\textbf{Acknowledgment:} The second author greatly acknowledges The Council of Scientific and Industrial Research (CSIR File No.: 09/025(0284)/2019-EMR-I), Government of India, for
the award of JRF.

                         \end{document}